\documentclass[12pt,a4paper]{article}
\usepackage{amssymb}
\usepackage{amsmath,amsfonts,amsthm}

\usepackage{graphicx}
\usepackage{amsmath}

\usepackage{amsfonts}
\usepackage{amsthm,amscd}
\usepackage{graphicx}
\usepackage{amsmath}
\usepackage{amssymb}
\usepackage{amstext}

\newtheorem{theorem}{Theorem}

\newtheorem{definition}[theorem]{Definition}

\newtheorem{lemma}[theorem]{Lemma}

\newtheorem{proposition}[theorem]{Proposition}

\title{Semiclassical limits, Lagrangian states and coboundary equations}

\author{Artur O. Lopes and Joana Mohr}

\begin{document}

\maketitle

\begin{abstract}

Assume that $f$ is  a continuous transformation $f:S^1 \to S^1$. We consider  here the cases where  $f$ is either the transformation $f(z)=z^2$ or $f$  is a smooth diffeomorphism of the circle $S^1$.

Consider a fixed continuous potential $\tau:S^1=[0,1)\to \mathbb{R}$,  $\nu\in \mathbb{R}$ and $\varphi:S^1 \to \mathbb{C}$ (a quantum state). The transformation
$\hat F_{\nu}$ acting on $\varphi:S^1 \to \mathbb{C}$, $\hat F_{\nu}(\varphi)=\phi$,  defined by
$\displaystyle \phi(z)=\hat F_{\nu}(\varphi(z))=\varphi(f(z))e^{i\nu\tau(z)}$
 describes a discrete time dynamical evolution of the quantum state $\varphi$.

Given
$S:\mathbb{R}\to \mathbb{R }$ we define the Lagrangian state
$$ \varphi_{x}^S(z)=\sum_{k\in\mathbb{Z}}e^{\frac{iS (z-k)}{\hbar}}e^{-\frac{(z-k-x)^2}{4\hbar}}.$$

 In this case $\hat F_{\nu}(\varphi_{x}^S(z))=\sum_{k\in\mathbb{Z}}e^{\frac{iS (f(z)-k)}{\hbar}}e^{-\frac{(f(z)-k-x)^2}{4\hbar}}e^{i\nu\tau(z)}$.

Under suitable conditions on $S$ the micro-support of $\varphi^S_x (z)$, when $\hbar \to 0$, is $(x,S'(x))$.
One  of meanings of the semiclassical limit in Quantum Mechanics  is to take  $\nu=\frac{1}{\hbar}$ and $\hbar \to 0$.
We address the question of finding $S$ such that $\varphi^S_x$ satisfies the property: $ \forall x$,  we have that $\hat{F}_\nu(\varphi^S_x)$ has micro-support on the graph of $y\to S'(y)$ (which is the micro-support of $\varphi^S_x$). In other words: which $S$ is such that $\hat{F}_\nu$ leaves the micro-support of $\varphi^S_x$ invariant? This is related to a coboundary equation for $\tau$,  twist conditions and   the boundary of the fat attractor.

\end{abstract}

\section{Introduction}

The Gaussian  wavepacket for $x\in \mathbb{R}, \xi\in \mathbb{R}$ is the function
$$y\to\tilde{\varphi}_{x,\xi}\, (y)=\,\, e^{ \, \frac{i\xi\,y}{\hbar}\, \,} e^{- \frac{(y-x)^2}{4\, a^2}}  .$$

This quantum state  describes a quantum particle located at $x$ and moment $\xi$.  The variance of the position is $a$ and the mean is $
x$ (see \cite{Sch} or \cite{Gia}).

We assume the variance of the position is $a=\sqrt{\hbar}$.
Then we get

\begin{equation}  \label{gp} \tilde{\varphi}_{x,\xi}\, (y)= \tilde{\varphi}_{x,\xi,\hbar}\, (y)= \,\, e^{ \, \frac{i\xi\,y}{\hbar}\, \,}\, \,e^{- \frac{(y-x)^2}{4\, \hbar}}.
\end{equation}

One can also consider the case where the variance of the position is $a=\sqrt{\hbar/m}$ and similar results for the setting we consider here
can also be obtained, when $m\to \infty$ and $\hbar$ is fixed. That is, $\hbar \to 0$ for the present setting is equivalent to $m\to \infty$ in this new setting.

In any case the Gaussian  wavepacket minimizes the Heisenberg uncertain principle for states with mean position $x$ and mean moment $\xi=m\, v$ (see \cite{Sch}).

When $\hbar \to 0$ the distribution of the position of  of $\tilde{\varphi}_{x,\xi}\, $ will be more and more concentrated in $x$. The  natural way to describe this property  which also contemplates the momentum is the concept of microsupport (see \cite{Mart} or \cite{Lop}).

As we are on a setting of functions which are defined on $S^1$ we will have to consider the periodic Gaussian wavepacket (see \cite{Faure}):
\begin{equation}  \label{gpp} \varphi_{x,\xi}(z)=\sum_{k\in\mathbb{Z}}\tilde \varphi_{x,\xi}(z-k)=\sum_{k\in\mathbb{Z}}e^{\frac{i\xi (z-k)}{\hbar}}e^{-\frac{(z-k-x)^2}{4\hbar}}.
\end{equation}



We say that a function $z(\hbar)$ is $O(\hbar^\infty)$ if  for each $ N>0$ there exists $C_N$ and $\delta>0$ such that
$$ |z(\hbar)|\leq C_N \, \hbar^N, \,\,\,\,\text{for}\,\,\,|\hbar| \leq \delta.$$

\medskip

\begin{definition} Given a  family of functions $\phi_\hbar $, $\hbar\sim 0$,  in $\mathcal{L}^2 (S^1)$, we say that it is micro-locally small near $(x_0,\xi_0)$ if
$$|< \,\varphi_{x,\xi}, \phi_\hbar>\, |$$
is $O(\hbar^{\infty})$ uniformly in a neighbourhood of  $(x_0,\xi_0)$. The complementary of such points  $(x_0,\xi_0)$ is called the micro-support of the family $\phi_\hbar $.

\end{definition}

Consider fixed $(x,\xi)$ and the corresponding periodic Gaussian wavepacket $\phi_\hbar=\varphi_{x,\xi}$ (according to (\ref{gpp})) depending on $\hbar$.
For any $(y,\eta)\neq (x,\xi) $ the modulus of the inner product
$$|< \varphi_{y,\eta}\,,\,\varphi_{x,\xi} > | =|< \varphi_{y,\eta}\,,\,\phi_\hbar > | $$
is $O(\hbar^\infty)$.

Therefore, the micro-support of  $\varphi_{x,\xi}\, $ is $(x,\xi)\, $. See the proof in chapter 3 in \cite{Mart} for the non-periodic case. The proof in the periodic case is similar (see also \cite{Faure}).

\bigskip

We consider a fixed continuous potential $\tau:S^1=[0,1)\to \mathbb{R}$ and $\nu\in \mathbb{R}$.
 We also consider a general dynamical system described by a continuous transformation $f:S^1 \to S^1$. One particular case is the transformation $f(z)=z^2$ (section 2). Another possibility we will consider here is when $f$ is a smooth diffeomorphism of the circle (section 3).

Then, for a given $z \in[0,1)$ and function $\varphi$ taking complex values we denote the transformation
$\hat F_{\nu}$, acting on $\varphi$,   by
$$\hat{F}_\nu(\varphi (z) )=\phi(z)= \varphi(f(z)) \, e^{\,i\,\nu\, \tau (z)}.$$

One can consider the above expression as a kind of discrete time dynamical evolution of a quantum state $\varphi$. The transformation $\hat{F}_\nu$ was  considered in \cite{Faure} and we borrow the notation of that paper.

Note that if $\varphi=\varphi_{x,\xi}$ is the Gaussian wavepacket, we have
$$\hat{F} (\varphi_{x,\xi}\,(z)\,)= \sum_{k\in\mathbb{Z}}e^{\frac{i\, \xi(f(z)-k)}{\hbar} }e^{-\frac{(f(z)-k-x)^2}{4\hbar}}e^{i\nu\tau(z)}.$$

\bigskip

Given $S:\mathbb{R} \to \mathbb{R}$ a Lagrangian state is an expression of the form
$$ a(x)\, e^{ \frac{i}{\hbar}\,S(x)},$$
where $a:\mathbb{R} \to \mathbb{R}$ is such that $\int a(x)^2 dx=1$.
One example is  $ \frac{1}{\sigma \sqrt{2 \pi}}\, e^{\frac{-x^2}{4 \sigma^2}}\, e^{ \frac{i}{\hbar}\,S(x)}.$

We consider here a periodic setting:
given
$S:S^1\to \mathbb{R }$ the associated  Lagrangian state is the function
$$ \varphi_{x}^S(z)=\sum_{k\in\mathbb{Z}}e^{\frac{iS (z-k)}{\hbar}}e^{-\frac{(z-k-x)^2}{4\hbar}}.$$

We are interested in the microsupport of  $\hat F_{\nu}(\varphi_{x}^S(z))$ in the semiclassical limit, $\nu=\frac{1}{\hbar}, \hbar \to 0$.
The analysis of the problem in the case of the transformation $f(z)=z^2$ is related to a certain coboundary equation. It is interesting that this question is associated to the study of the boundary of a certain attractor.

The general study of the boundary of attractors is the object of several papers in dynamics (see \cite{Nu} and \cite{LO})

We would like to thanks D. Smania for helpful conversations related to section \ref{secdif} of the present paper.
We would like to thanks J. Hounie for providing us with some nice references.

\section{The transformation $f(x)=2 \,x $ (mod 1)}

\bigskip

Consider a fixed continuous function $A: S^1=[0,1) \to \mathbb{R}$.
Remember that $f(z)=z^2$ (or, in an equivalent way $f(x)=2\, x$ (mod 1)).

Consider $F :S^1 \times \mathbb{R} \to S^1 \times \mathbb{R}$ given by
$$F(z,s)= (z^2, \lambda \,s + A(z))= (f(z), \lambda\, s+ A(z)),$$
 where $0<\lambda\leq 1$. The term $\lambda$ is called the discounted term. In \cite{LO} the ergodic analysis of this problem when $\lambda$ is fixed was considered and was also analyzed the  limit problem when $\lambda \to 1$. Here we consider the case when $\lambda=\frac{1}{2}$.

 We denote by $G_1$ and $G_2$ the two inverse branches of $F$, which are respectively
 $$G_1(y,r)= \bigg(\frac{y}{2}, \frac{r- A(\frac{y}{2})}{\lambda }\bigg) \,\,\,\text{and}\,\,\,G_2(y,r)= \bigg(\frac{y}{2}\,\,+\,\frac{1}{2}, \frac{r- A(\frac{y}{2}\,\,+\,\frac{1}{2})}{\lambda }\bigg) $$

Among others possible choices are $\lambda= \frac{1}{2}$ (preserves volume), or  $\lambda=1.$

In Proposition 6 in  \cite{Faure} the case when $A=0$ and $\lambda=\frac{1}{2}$ is considered.

\begin{theorem} (see \cite{Faure}) Suppose $f(z)=2\, z$ (mod 1). As $\hbar \to 0$,
$$\phi_\hbar(z)=\hat{F}_\nu(\varphi_{x,\xi}\, (z))= \varphi_{x,\xi}\, (f(z))\, e^{i\,  \nu \tau(z)}$$

has micro-support in
$$ \{ G_1(x,\xi), G_2 (x,\xi)\}\subset T^*\, S^1,$$
where $G_1$ and $G_2$ are the inverse branches  for
$$F(z,s)= (z^2, \lambda \,s)= (f(z), \frac{1}{2}\, s ).$$

\end{theorem}

{\bf Proof:}

For fixed $(x,\xi)$ and variable $(y,\eta)$, $x,y\in S^1$
$$< \varphi_{y,\eta}\,,\,\hat{F}_\nu\varphi_{x,\xi} >= < \varphi_{y,\eta}\,,\varphi_{x,\xi}\, (f)\, e^{i\,  \nu \tau}>= $$

$$ \int_{S^1}\,\overline{ \varphi_{y,\eta}\, (z)}\,\sum_{k\in\mathbb{Z}}e^{\frac{i\, \xi(f(z)-k)}{\hbar} }e^{-\frac{(f(z)-k-x)^2}{4\hbar}}e^{i\nu\tau(z)}\,dz= $$
$$ \int_{S^1}\,\bigg(\sum_{q\in\mathbb{Z}}e^{-\frac{i\eta (z-q)}{\hbar}}e^{-\frac{(z-q-y)^2}{4\hbar}}\,\bigg)\,\bigg(\sum_{k\in\mathbb{Z}}e^{\frac{i\, \xi(f(z)-k)}{\hbar} }e^{-\frac{(f(z)-k-x)^2}{4\hbar}}e^{i\nu\tau(z)}\bigg)\,dz. $$
Now if we define $S(z)=\xi z$,  then we can apply proposition \ref{micros} (see appendix), with a slight modification, because $e^{i\nu\tau(z)}$ does not depend on $\hbar$, to obtain that
$$< \varphi_{y,\eta}\,,\,\hat{F}_\nu\varphi_{x,\xi} >=O(\hbar^{\infty})$$
 if $f(y)\neq x$ or if $f(y)= x$ and $\eta\neq \xi f'(y)=2\xi$.
 Then the micro-support requires $ f(y)=x$ and  $\eta=2\, \xi$.  That is,       $(y,\eta)$ is contained in
$\{ G_1(x,\xi), G_2 (x,\xi)\}$.
\qed

\bigskip

Now we assume that $\nu $ and $\hbar$ are coupled by $\hbar=\frac{1}{\nu}$
and $A=-\frac{1}{2}\tau'$.
\medskip

\begin{theorem} (see \cite{Faure})  Suppose $f(z)=2\, z$ (mod 1).  As $\hbar=\frac{1}{\nu}\,\to 0$
$$\phi_\hbar (z) =\hat{F}_\nu(\varphi_{x,\xi}\, (z))= \varphi_{x,\xi}\, (f(z))\, e^{i\,  \nu \tau(z)} =\varphi_{x,\xi}\, (f(z))\, e^{  \frac{i\,\tau(z)}{\hbar}}$$

has micro-support in
$$ \{ G_1(x,\xi), G_2 (x,\xi)\}\subset T^*\, S^1,$$
where $G_1$ and $G_2$ are the inverse branches  for
$$F(z,s)= (z^2, \frac{1}{2}\, s - \frac{1}{2}\tau'(z)) = (f(z), \frac{1}{2}\, s - \frac{1}{2}\tau'(z)).$$

\end{theorem}

{\bf Proof:}

For fixed $(x,\xi)$ and variable $(y,\eta)$, $x,y\in S^1$
$$< \varphi_{y,\eta}\,,\,\hat{F}_\nu\varphi_{x,\xi} >= < \varphi_{y,\eta}\,,\varphi_{x,\xi}\, (f)\, e^{i\,  \nu \tau}>= $$
$$= \int_{S^1}\,\bigg(\sum_{q\in\mathbb{Z}}e^{-\frac{i\eta (z-q)}{\hbar}}e^{-\frac{(z-q-y)^2}{4\hbar}}\,\bigg)\,\bigg(\sum_{k\in\mathbb{Z}}e^{\frac{i\, \xi(f(z)-k)}{\hbar} }e^{-\frac{(f(z)-k-x)^2}{4\hbar}}e^{\frac{i\tau(z)}{\hbar} }\bigg)\,dz. $$

Again, if we define $S(z)=\xi z$,  then we can apply proposition \ref{micros} (see appendix) to obtain that
$$< \varphi_{y,\eta}\,,\,\hat{F}_\nu\varphi_{x,\xi} >=O(\hbar^{\infty})$$
 if $f(y)\neq x$ or if $f(y)= x$ and $\eta\neq \xi f'(y)+\tau'(y)=2\xi+\tau'(y)$.
 Then the micro-support requires $ f(y)=x$ and  $\eta=2\, \xi+\tau'(y)$.  That is,       $(y,\eta)$ is contained in
$\{ G_1(x,\xi), G_2 (x,\xi)\}$.
\qed

\bigskip

Now we consider a different kind of problem.

If
$S:\mathbb{R }\to \mathbb{R }$ we define the Lagrangian state  as the periodic function

 $$ \varphi_{x}^S(z)=\sum_{k\in\mathbb{Z}}e^{\frac{iS (z-k)}{\hbar}}e^{-\frac{(z-k-x)^2}{4\hbar}}.$$


\bigskip
It is a well known fact that:

\begin{proposition}
The Lagragian state $\varphi_{x}^S $, associated to $S$, has microsupport on  $(x,S'(x))$, when $\hbar  \to 0.$

\end{proposition}

{\bf Proof}:
For fixed $x$ and variable $y,\eta$ we have

$$< \varphi_{y,\eta}  ,\varphi_{x}^S>=\int_{S^1} \overline{\sum_{q\in\mathbb{Z}}e^{\frac{i\eta (z-q)}{\hbar}}e^{-\frac{(z-q-y)^2}{4\hbar}}  }\,\sum_{k\in\mathbb{Z}}e^{\frac{iS (z-k)}{\hbar}}e^{-\frac{(z-k-x)^2}{4\hbar}}\, dz.$$

It follows by the same kind of  arguments used in Proposition \ref{micros}, that the micro-support requires $y=x$ and $\eta=S'(x)$, i.e.,
$\varphi_{x}^S$
has microsupport on the graph $(x,S'(x))$, when $\hbar\to 0$, as we claimed.

\qed

\bigskip

Lagrangian states of the form
$u_\hbar(x)= a(x) e^{\frac{iS(x)}{\hbar} \, }$ were consider in Example 2.3 in page 32 in \cite{Anan}, where $a$ is differentiable and positive, in the setting of Aubry-Mather Theory.
\medskip

 Now, if we consider $f:S^1\to S^1$ we define
 $$\hat F_{\nu}(\varphi_{x}^S(z))=\varphi_{x}^S(f(z))e^{i\nu\tau(z)}=\sum_{k\in\mathbb{Z}}e^{\frac{iS (f(z)-k)}{\hbar}}e^{-\frac{(f(z)-k-x)^2}{4\hbar}}e^{i\nu\tau(z)}.$$

\begin{definition}
Given a continuous function $A:S^1\to \mathbb{R}$ and $\lambda\in (0,1)$, we say that a continuous function
{\bf $b=b_\lambda:[0,1) \to \mathbb{R}$ is a $\lambda$-calibrated plus-subaction for $A$}, if for all $z\in S^1$,
$ b(z) = \displaystyle\max_{f(y)=z} \{ \lambda \, b(y) + A(y)\},$ where $f(z)=2\, z$ (mod 1).
\end{definition}

This problem is related to $\lambda$-maximizing probabilities (see \cite{LO}).

 We will consider here the case where $\lambda=\frac 1 2$ and $A=-\frac{\tau'}{2}.$

Denote by $u$ the $\frac 1 2$-subaction for $A(z)=-\frac{1}{2}\tau'(z)$. The graph $(z, u(z))$ is the upper boundary of the  fat attractor and invariant for  $F$, where
$$F(z,s)= \Big(z^2, \frac{1}{2}\, \,s- \frac{1}{2}\tau'(z)\Big)= \Big(f(z), \frac{1}{2}\, s - \frac{1}{2}\tau'(z)\Big).$$
Then we have $(y,u(y))=F(z,u(z))=\Big(f(z), \frac{1}{2}\, u(z) - \frac{1}{2}\tau'(z)\Big)$, this implies $y=f(z)$ and $\frac{1}{2} u(z) - \frac{1}{2}\tau'(z)=u(y)=u(f(z))$,
in this way for any $z$

\begin{equation} \label{se} u(f(z))=  \frac{1}{2} u(z) - \frac{1}{2} \tau '(z).
\end{equation}

We call the equation $\frac{1}{2} \tau '(z)\,=\,   \frac{1}{2} u(z) - u(f(z))$ a {\bf $\frac{1}{2}$-coboundary equation} for $u$ and $\tau'$.

In several examples $u$ is $C^\infty$ up to some finite number of points (this kind of property is the main issue in \cite{LO}). On these points there exists the left and right limit of the derivative (see \cite{LO}) which are not zero nor infinite. This happens, for instance, if the potential $A$ is $C^\infty$ and $A=\tau'$ satisfies the twist condition to be defined later (see \cite{LO}).

In the appendix in
the estimates of the asymptotic $\hbar\to 0$ it will be important for technical reasons that $u$ is piecewise smooth $C^\infty$ (left and right finite non zero derivatives in a finite number of singular points).


From now on we assume $\tau$ is $C^\infty$ and then can write
locally $u = S'$ for some $S:[0,1) \to \mathbb{R}$ which is $C^\infty$ by parts.

Then, from (\ref{se}) we get that for any
$z$:
$$ S'(f(z))=  \frac{1}{2} S'(z) - \frac{1}{2} \tau '(z),$$

or
\begin{equation} \label{seS} 2\, S'(f(z)) =  S'(z) - \tau '(z).\end{equation}

\bigskip

We would like to consider different possible $S:[0,1) \to \mathbb{R}$, compare the corresponding $S$-Lagrangian state and look for the $S$ such the associated micro-support (the graph of $S'$)   is  invariant under the action of $\hat{F}_\nu$, when $\hbar=\frac{1}{\nu} \to 0.$ Here  $\hat{F}_\nu$ is such that $\hat{F}_\nu(\phi\, (z))= \phi\, (f(z))\, e^{i\,  \nu \tau(z)}$.

The $S$ such that  $ S'=u$, where $u$ satisfies (\ref{se}), is the one we have to look for as we will show now.

\medskip

Remember that for any $x$ we have that $\varphi_{x}^S$ has micro-support on the graph $(x,S'(x))$.

\medskip

\begin{theorem} \label{subsup}

Consider
$$\,\varphi_{x}^S\, (z)= \sum_{k\in\mathbb{Z}}e^{\frac{iS (z-k)}{\hbar}}e^{-\frac{(z-k-x)^2}{4\hbar}}$$
where $S'=u$ is the $\frac{1}{2}$-plus-calibrated subaction  for $-\frac{\tau'}{2}$ (see equation \ref{se}).
Assume that $S$ is piecewise smooth $C^\infty$ (left and right finite non zero derivatives in a finite number of singular points), and also that $\tau$ is $C^{\infty}$.

When $\hbar=1/\nu\,\to 0$, we get that
$$\hat{F}_\nu(\varphi_{x}^S\, (z))= \varphi_{x}^S\, (f(z))\, e^{i\,  \nu \tau(z)} =\varphi_{x}^S\, (f(z))\, e^{\,  \frac{i\tau(z)}{\hbar}\,}$$
has  micro-support on the graph of
$S'$ (locally defined).

The underlying dynamics is
$$F(z,s)= \Big(f(z), \frac{1}{2}\, s - \frac{1}{2}\tau'(z)\Big).$$

\end{theorem}

{\bf Proof:}
For fixed $x$ and variable $\eta,y$, where $x,y\in S^1$, first we apply proposition \ref{micros} (see appendix) to obtain that
$$< \varphi_{y,\eta}\,,\,\hat{F}_\nu\varphi_{x}^S >=O(\hbar^{\infty})$$
 if $f(y)\neq x$ or if $f(y)= x$ and $\eta\neq S'(f(y)) f'(y)+\tau'(y)=2S'(f(y))+\tau'(y)$.
Then the micro-support requires $ f(y)=x$ and $\eta=2\, S'(f(y)) + \tau'(y)$. As we assume that $S'$ is $\frac{1}{2}-$subaction for
$  - \frac{1}{2} \tau ' $, then,  from (\ref{seS}) we get that $\eta = S'(y)=u(y)$. Therefore micro-support is   on the graph  $(y,u(y))=(f^{-1}(x),u(f^{-1}(x)))$ of $u=S'$.


\qed

\bigskip

Now we will present sufficient conditions on $\tau$ in order the function $S$ is under the hypothesis of last theorem.

\bigskip

Denote $A=- \frac{1}{2} \tau'$. The inverses branches of $f $ are denoted by $\tau_1$ and $\tau_2$.
We assume that $\tau_1(0,1)=(0,0.5)$ and  $\tau_2(0,1)=(0.5,1)$.

 We
also denote $a=(a_1,a_2,...)$ a generic element in $ \Omega = \{1,2\}^\mathbb{N}$.

Consider (as Tsujii in \cite{T1}) the function $s: (S^1, \Omega) \to \mathbb{R}$, where $ \Omega = \{1,2\}^\mathbb{N}$, given by
$$ s(x,a) = \sum_{k=0}^{\infty} \Big(\frac{1}{2}\Big)^k A(\,(\tau_{a_k}\circ \tau_{a_{k-1}}\circ\,...\, \circ\tau_{a_0})\, (x)  \,) ,$$
and, $a=(a_0,a_1,a_2,...).$

Note that for a fixed $a$ the function $s(x,a)$ is not a  periodic function on $x\in[0,1]$.

We define $\pi(x)=i$, if $x$ is in the image of
$\tau_i(S^1)$, $i\in \{1,2\}$.

Note that $(x,a) \to (f(x),  \pi(x) a )$ defines a skew product on $S^1\times \Omega$.

Note also \cite{T1} that
$$s(f(x), \pi(x) a)= A(x) + \frac{1}{2}\, s(x,a).$$

\bigskip

For a fixed $b\in \Omega$ consider
$$\varphi_{x}^{S(\cdot,b)}\, (z)=  \,\, e^{ \, \frac{iS(z,b)}{\hbar}\, \,\,}\, \,e^{- \frac{(z-x)^2}{4\, \hbar}}.$$
If we suppose $s(z,b)=\frac{\partial S}{\partial z}(z,b)$, then we can show in a similar way as in last result that
$$\hat{F}_\nu(\varphi_{x}^{S(\cdot,b)}\, (z))$$
has support on the graph $(\,y, s(y ,a)\,)$ where $f(y)=x$ and $\pi(y)\,a=b.$

\bigskip

From \cite{LO} one gets that
the plus-calibrated subaction $u$ satisfies
$$u(x)= \sup_{ c \in \{1,2\}^\mathbb{N}} \, s(x,c).$$
This property
resembles the classical method of  obtaining solutions of the Hamilton-Jacobi equation via envelopes
(Huygens's generation of the wave front) as described for instance  in \cite{FJ} Chapter I.9.

\bigskip

\begin{definition}
Consider a fixed $\overline{x}\in S^1$ and variable $x\in S^1$, $a\in\{1,2\}^\mathbb{N}$, then we
define
\begin{equation} \label{W}
W(x,a)= s(x,a) - s(\overline{x},a).
\end{equation}
We call such $W$ the $\frac{1}{2}$-involution kernel for $A=-\frac{1}{2}\tau'$.
\end{definition}


For $a$ fixed $W(x,a)$ is smooth on $x\in(0,1).$

Below we consider the lexicographic order in $\{1,2\}^\mathbb{N}$.
\begin{definition} \label{tui} We say that $A$ satisfies  the twist condition, if an (then, any) associated involution kernel $W$, satisfies the property: for any $a<b$, we have
$$\frac{\partial W}{\partial x}(x,a) - \frac{\partial W}{\partial x}(x,b)>0.$$
\end{definition}

It is equivalent to state the above relation for $s$ or for $W$.

If $A=-\frac{1}{2}\tau'$ satisfies the twist condition, then, $u$ will be piecewise smooth (see Corollary 15 in \cite{LO} ) and Theorem \ref{subsup} can be applied.

\bigskip

\bigskip

\section{Diffemorphisms of the circle $f:S^1 \to S^1$} \label{secdif}

\bigskip

Suppose $f$ is an orientation preserving diffeomorphism of the circle $S^1$ of class $C^k$ and $\tau:S^1 \to \mathbb{R}$ is also of class $C^k$. The value of $k$ may depend of $f$ but will require at least $k\geq 4.$

We are interest in condition on $f$ and $\tau$ such that there exists
 $u:S^1 \to \mathbb{R}$ which is differentiable and  satisfies
\begin{equation} \label{equacao u}
-\tau'(z)=u(f(z))f'(z)-u(z).\end{equation}



Suppose $f$ as above and $A:S^1 \to \mathbb{R}$ is also of class $C^k$. If there exists differentiable  $w:S^1 \to \mathbb{R}$ such that
\begin{equation} \label{equacao A}A(z)=w(f(z))-w(z),\end{equation}
then, taking derivative on both sides, denoting $u=w'$ and $\tau=-A$, we get (\ref{equacao u}).

Given a smooth $A$, the problem of finding a smooth $w$ is considered in section 2.2 in \cite{AK} and \cite{Her}.

Given a real number
$\alpha$  we say that it satisfies the Diophantine condition if: $$\bigg|\alpha-\frac{p}{q}\bigg|>\frac{K}{q^{2+\beta}} \mbox {   for all  } \frac{p}{q}\in \mathbb{Q},$$ where $K$ and $\beta$ are positive constants.

In the case $f(z)=R_\alpha(z)=z+\lambda$, where $\lambda$ is diophantine and $A=-\tau$ is $C^{3+\alpha}$, there  exists a differentiable function $w$
such that
\begin{equation}\label{eqrot}
-\tau (z)=w(f(z))-w(z),
\end{equation}
as before,  taking $u=w'$  we get a solution for equation (\ref{equacao u}).

\bigskip

Now let $f$ be a $ C^k$ circle diffeomorphism, $k\geq 3$ and suppose the rotation number $\alpha$ of $f$ satisfies the Diophantine condition: $$\bigg|\alpha-\frac{p}{q}\bigg|>\frac{K}{q^{2+\beta}} \mbox {   for all  } \frac{p}{q}\in \mathbb{Q},$$ where $K$ and $\beta$ are positive constants. Suppose also that $k>2\beta+1$, then, given $A=-\tau$ of class $C^k$ there exists $w$ of class $C^k$ satisfying equation (\ref{equacao A}),
 again if $u=w'$  we get a solution for equation (\ref{equacao u}).

\bigskip


\begin{theorem} Suppose $f$ is a $C^\infty$ diffeomorphism of the circle,
 $\tau:S^1 \to \mathbb{R}$ is smooth $C^\infty$ and assume that for any $z$ the coboundary equation
\begin{equation}\label{equacao u1}
u(z)=u(f(z))\,f'(z)+\tau'(z)
\end{equation}
is true ($u$ is $C^\infty$).

Suppose also that $S'=u$.
Then, the micro-support of the Lagrangian state  $\varphi_{x}^S(y)=e^{\frac{iS( y)}{\hbar}}e^{-\frac{(y-x)^2}{4\hbar}}$ is invariant for the transformation $\hat{F_{\nu}}(\varphi_{x}^S(z))=\varphi_{x}^S(f(z))e^{\frac{\, i \tau (z)}{\hbar}},$
when $\hbar \to 0$.

The underlying dynamics is given by  $F(z,s)=\Big(f(z),\frac{s-\tau'(z)}{f'(z)}\Big)$,  $F:S^1 \times \mathbb{R} \to S^1 \times \mathbb{R}$ and the graph of $u$ is invariant by $F$.

\end{theorem}

{\bf Proof:}

As before, for fixed $x$ and variable $\eta,y$, where $x,y\in S^1$, first we apply proposition \ref{micros} (see appendix) to obtain that
$$< \varphi_{y,\eta}\,,\,\hat{F}_\nu\varphi_{x}^S >=O(\hbar^{\infty})$$
 if $f(y)\neq x$ or if $f(y)= x$ and $\eta\neq S'(f(y)) f'(y)+\tau'(y)$.

Then the micro-support requires $ f(y)=x$ and $\eta= S'(f(y))  f'(y) + \tau'(y)$.

\noindent And, as we suppose   that $S'=u$,  using equation \eqref{equacao u}, we obtain $$\eta=S'(f(y)) f'(y) +\tau'(y)=u'(f(y)) f'(y) +\tau'(y)=u(y),$$

Therefore, we get that the microsupport of $\hat{F_{\nu}}(\varphi_{x}^S)$ is on the graph  $(y,u(y))=(f^{-1}(x),u(f^{-1}(x)))$.

Finally, let us show that the graph of $u$ is invariant for $F$. In fact, note that if $y=f(z)$
then by \eqref{equacao u} we get $u(y)=u(f(z))=\frac{u(z)-\tau'(z)}{f'(z)}$ and
$$F(z,u(z))=\bigg(f(z),\frac{u(z)-\tau'(z)}{f'(z)}\bigg)=(y,u(y)). $$

\qed

\bigskip


\section{Appendix}

Let $(x,\xi)\in\mathbb{R}^2$ , then for each $z\in\mathbb{R}$ we define before  $$\tilde \varphi_{x,\xi}(z)=e^{\frac{i\xi z}{\hbar}}e^{-\frac{(z-x)^2}{4\hbar}},$$
and if
$z\in S^1$ we define the periodic function

 $$ \varphi_{x,\xi}(z)=\sum_{k\in\mathbb{Z}}\tilde \varphi_{x,\xi}(z-k)=\sum_{k\in\mathbb{Z}}e^{\frac{i\xi (z-k)}{\hbar}}e^{-\frac{(z-k-x)^2}{4\hbar}}.$$

 For each $x\in \mathbb{R}$ , $S:\mathbb{R}\to \mathbb{R }$ the Lagrangian package is defined,  for each  $z\in\mathbb{R}$, by
   $$\tilde \varphi_{x}^S(z)=e^{\frac{iS(z)}{\hbar}}e^{-\frac{(z-x)^2}{4\hbar}}.$$

   Remember that for
$S:\mathbb{R }\to \mathbb{R }$ we define the periodic function

 $$ \varphi_{x}^S(z)=\sum_{k\in\mathbb{Z}}\tilde \varphi_{x}^S(z-k)=\sum_{k\in\mathbb{Z}}e^{\frac{iS (z-k)}{\hbar}}e^{-\frac{(z-k-x)^2}{4\hbar}}.$$

 And, if we take $f:S^1\to S^1$ we define
 $$\hat F_{\nu}(\varphi_{x}^S(z))=\varphi_{x}^S(f(z))e^{i\nu\tau(z)}=\sum_{k\in\mathbb{Z}}e^{\frac{iS (f(z)-k)}{\hbar}}e^{-\frac{(f(z)-k-x)^2}{4\hbar}}e^{i\nu\tau(z)}$$

Therefore, if  $z\in S^1$

 $$<\varphi_{y,\eta}(z),\hat F_{\nu}(\varphi_{x}^S(z)) >=$$
 $$\int_{S^1} \sum_{j\in\mathbb{Z}}e^{\frac{i\eta (z-j)}{\hbar}}e^{-\frac{(z-j-y)^2}{4\hbar}} \sum_{k\in\mathbb{Z}}e^{-\frac{iS (f(z)-k)}{\hbar}}e^{-\frac{(f(z)-k-x)^2}{4\hbar}}e^{i\nu\tau(z)}dz $$

\begin{proposition}\label{estim} Let us fix $x,y\in \mathbb{R}$. Then, there exist constants $\bar k $ and $\bar j$ such that
$$<\varphi_{y,\eta}(z),\hat F_{\nu}(\varphi_{x}^S(z)) >\approx\int_{S^1} e^{\frac{i\eta (z+\bar j)}{\hbar}}e^{-\frac{(z+\bar j-y)^2}{4\hbar}} e^{-\frac{iS (f(z)+\bar k)}{\hbar}}e^{-\frac{(f(z)+\bar k-x)^2}{4\hbar}}e^{i\nu\tau(z)}dz,  $$
\end{proposition}

{\bf Proof:}
We  choose the constants  $\bar k, \bar l$ such that  $x\in(\bar k,\bar k+1 )$ and $y\in(\bar j ,\bar j+1 )$.

We will show that the other terms in the integral rapidly  decrease  to zero, when   $\hbar\to 0$. Indeed,
    $$<\varphi_{y,\eta}(z),\hat F_{\nu}(\varphi_{x}^S(z)) >=$$
    $$=\int_{S^1} \sum_{j\in\mathbb{Z}}e^{\frac{i\eta (z-j)}{\hbar}}e^{-\frac{(z-j-y)^2}{4\hbar}} \sum_{k\in\mathbb{Z}}e^{-\frac{iS (f(z)-k)}{\hbar}}e^{-\frac{(f(z)-k-x)^2}{4\hbar}}e^{i\nu\tau(z)}dz =$$
  $$=\int_{S^1} e^{\frac{i\eta (z+\bar j)}{\hbar}}e^{-\frac{(z+\bar j-y)^2}{4\hbar}} e^{-\frac{iS(f(z)+\bar k)}{\hbar}}e^{-\frac{(f(z)+\bar k-x)^2}{4\hbar}}e^{i\nu\tau(z)}dz +$$
    $$ + \int_{S^1} e^{-\frac{iS (f(z)+\bar k)}{\hbar}}e^{-\frac{(f(z)+\bar k-x)^2}{4\hbar}}e^{i\nu\tau(z)}\sum_{j\in\mathbb{Z}, j\neq -\bar j}e^{\frac{i\eta (z-j)}{\hbar}}e^{-\frac{(z-j-y)^2}{4\hbar}} dz +$$
    $$ + \int_{S^1} e^{\frac{i\eta (z+\bar j)}{\hbar}}e^{-\frac{(z+\bar j-y)^2}{4\hbar}} e^{i\nu\tau(z)}\sum_{k\in\mathbb{Z},k\neq -\bar k}e^{-\frac{iS (f(z)-k)}{\hbar}}e^{-\frac{(f(z)-k-x)^2}{4\hbar}}dz +$$
$$+\int_{S^1} \sum_{j\in\mathbb{Z}, j\neq -\bar j}e^{\frac{i\eta (z-j)}{\hbar}}e^{-\frac{(z-j-y)^2}{4\hbar}} \sum_{k\in\mathbb{Z},k\neq -\bar k}e^{-\frac{iS (f(z)-k)}{\hbar}}e^{-\frac{(f(z)-k-x)^2}{4\hbar}}e^{i\nu\tau(z)}dz .$$

We define the positive constant $$C=\min\{d(y,\bar j),d(y,\bar j+1), d(x,\bar k),d(x,\bar k+1)\}.$$ As  $z\in S^1$, we obtain
$$ \sum_{j\in\mathbb{Z}, j\neq -\bar j}e^{-\frac{(z-j-y)^2}{4\hbar}}=e^{-\frac{(z+\bar j-1-y)^2}{4\hbar}}+e^{-\frac{(z+\bar j+1-y)^2}{4\hbar}}+e^{-\frac{(z+\bar j-2-y)^2}{4\hbar}}+e^{-\frac{(z+\bar j+2-y)^2}{4\hbar}}+... $$
$$\leq 2( e^{-\frac{C^2}{4\hbar}}+e^{-\frac{(C+1)^2}{4\hbar}}+e^{-\frac{(C+2)^2}{4\hbar}}+...)\leq 2\bigg(e^{-\frac{C^2}{4\hbar}}+\int_C^{\infty}e^{-\frac{x^2}{4\hbar}}dx\bigg)=O\Big(e^{-\frac{C^2}{4\hbar}}\Big).$$

Therefore, the second term above can be bounded by

$$\bigg|\int_{S^1} e^{-\frac{iS (f(z)+\bar k)}{\hbar}}e^{-\frac{(f(z)+\bar k-x)^2}{4\hbar}}e^{i\nu\tau(z)}\sum_{j\in\mathbb{Z}, j\neq -\bar j}e^{\frac{i\eta (z-j)}{\hbar}}e^{-\frac{(z-j-y)^2}{4\hbar}} dz\bigg|\leq$$
$$\leq \int_{S^1} \sum_{j\in\mathbb{Z}, j\neq -\bar j}e^{-\frac{(z-j-y)^2}{4\hbar}} dz=O\Big(e^{-\frac{C^2}{4\hbar}}\Big)$$

The third term above can be bounded by the same arguments.

Note that if  $\hbar$ small enough we have that
$\displaystyle \sum_{j\in\mathbb{Z}, j\neq -\bar j}e^{-\frac{(z-j-y)^2}{4\hbar}}\leq 1$.
Therefore, the last term above can be bounded by

$$\bigg |\int_{S^1} \sum_{j\in\mathbb{Z}, j\neq -\bar j}e^{\frac{i\eta (z-j)}{\hbar}}e^{-\frac{(z-j-y)^2}{4\hbar}} \sum_{k\in\mathbb{Z},k\neq -\bar k}e^{-\frac{iS (f(z)-k)}{\hbar}}e^{-\frac{(f(z)-k-x)^2}{4\hbar}}e^{i\nu\tau(z)}dz \bigg |\leq$$
$$\leq \int_{S^1} \sum_{j\in\mathbb{Z}, j\neq -\bar j}e^{-\frac{(z-j-y)^2}{4\hbar}} \sum_{k\in\mathbb{Z},k\neq -\bar k}e^{-\frac{(f(z)-k-x)^2}{4\hbar}}dz\leq$$
$$\leq \int_{S^1}  \sum_{k\in\mathbb{Z},k\neq -\bar k}e^{-\frac{(f(z)-k-x)^2}{4\hbar}}dz=O\Big(e^{-\frac{C^2}{4\hbar}}\Big).$$

 \qed

Now, if  $x,y\in S^1$ and $\nu=\frac{1}{\hbar}$, then we get 

 $$<\varphi_{y,\eta}(z),\hat F_{\nu}(\varphi_{x}^S(z)) >\approx\int_{S^1} e^{-\frac{i\eta z}{\hbar}}e^{-\frac{(z-y)^2}{4\hbar}} e^{\frac{iS (f(z))}{\hbar}}e^{-\frac{(f(z)-x)^2}{4\hbar}}e^{\frac{i\tau(z)}{\hbar}}dz.  $$

\bigskip

\begin{lemma}\label{estderiv} Let us fix a constant $a$, then
$$\int_{y-a}^{y+a}\bigg|\frac{d^{2n}}{dz^{2n}}e^{-\frac{(z-y)^2}{4\hbar}}\bigg|dz$$ decay fast to zero as $\hbar\to 0$, for each $n\geq 1$.

Moreover, for each $n\geq 0$, we get
$$\int_{y-a}^{y+a}\bigg|\frac{d^{2n+1}}{dz^{2n+1}}e^{-\frac{(z-y)^2}{4\hbar}}\bigg|dz=O\Big(\frac{1}{\hbar^{n}}\Big).$$
\end{lemma}

\begin{proposition}\label{micros} Assume that $S$ is piecewise smooth $C^\infty$ (left and right finite non zero derivatives in a finite number of singular points).
Suppose $f$ is  $C^\infty$, $\nu=\frac{1}{\hbar}$ and
 $\tau:S^1 \to \mathbb{R}$ is smooth $C^\infty$. Let us fix $x,y\in S^1$. If $f(y)\neq x $ or $f(y)=x$ and $\eta\neq S'(f(y))f'(y)+\tau'(y)$, then $$<\varphi_{y,\eta}(z),\hat F_{\nu}(\varphi_{x}^S(z)) >=O(\hbar^{\infty}).$$

\end{proposition}

{\bf Proof:}
Define $\phi(z)=S(f(z))-\eta z+\tau(z)$, then by the Proposition \ref{estim} it remains to  analyze the following  integral.

$$\int_{S^1}e^{-\frac{(z-y)^2+(f(z)-x)^2}{4\hbar}}e^{\frac{i(S(f(z))-\eta z+\tau(z))}{\hbar}}dz=\int_{S^1}e^{-\frac{(z-y)^2+(f(z)-x)^2}{4\hbar}}e^{\frac{i\phi(z)}{\hbar}}dz$$

\bigskip
\textbf{Case 1:} We suppose that $f(y)\neq x$, then we can obtain  constants $C$ and $\epsilon$ such that
$(f(z)- x)^2>C^2$ for all $z\in (y-\epsilon, y+\epsilon)$, also if $z\in (y+\epsilon,1)$ or $z\in (0,y-\epsilon)$ we have that $(z-y)^2\geq \epsilon^2$.
Therefore, if we denote by $\delta=\min\{C,\epsilon\}$ we get
$$\bigg|\int_{S^1}e^{-\frac{(z-y)^2+(f(z)-x)^2}{4\hbar}}e^{\frac{i\phi(z)}{\hbar}}dz\bigg|\leq$$
$$\int_{S^1}e^{-\frac{(z-y)^2}{4\hbar}}e^{-\frac{(f(z)-x)^2}{4\hbar}}dz\leq e^{-\frac{\epsilon^2}{4\hbar}}\int_{0}^{y-\epsilon}e^{-\frac{(f(z)-x)^2}{4\hbar}}dz+$$
$$+e^{-\frac{C^2}{4\hbar}}\int_{y-\epsilon}^{y+\epsilon}e^{-\frac{(z-y)^2}{4\hbar}}dz+e^{-\frac{\epsilon^2}{4\hbar}}\int_{y+\epsilon}^{1}e^{-\frac{(f(z)-x)^2}{4\hbar}}dz=O\Big(e^{-\frac{\delta^2}{4\hbar}}\Big).$$

\bigskip
\textbf{Case 2:}  Suppose  that $f(y)=x$ and $\eta\neq S'(f(y))f'(y)+\tau'(y)$, i.e., $\phi'(y)=C\neq 0$, then  we can find a constant  $a$ such that $ \phi'(z)\geq \frac{C}{2}$ for  $z\in(y-a,y+a)$ Lebesgue a.e.w. Hence,

$$\int_{S^1}e^{-\frac{(z-y)^2+(f(z)-x)^2}{4\hbar}}e^{\frac{i\phi(z)}{\hbar}}dz=\int_{S^1}e^{-\frac{(z-y)^2}{4\hbar}}e^{\frac{i\phi(z)}{\hbar}}dz=$$
$$= \int_{0}^{y-a}e^{-\frac{(z-y)^2}{4\hbar}}e^{\frac{i\phi(z)}{\hbar}}dz+\int_{y-a}^{y+a}e^{-\frac{(z-y)^2}{4\hbar}}e^{\frac{i\phi(z)}{\hbar}}dz+\int_{y+a}^{1}e^{-\frac{(z-y)^2}{4\hbar}}e^{\frac{i\phi(z)}{\hbar}}dz .$$
Note that,  $\displaystyle\int_{0}^{y-a}e^{-\frac{(z-y)^2}{4\hbar}}e^{\frac{i\phi(z)}{\hbar}}dz=O\Big(e^{-\frac{a^2}{4\hbar}}\Big)$, $\displaystyle\int_{y+a}^{1}e^{-\frac{(z-y)^2}{4\hbar}}e^{\frac{i\phi(z)}{\hbar}}dz=O\Big(e^{-\frac{a^2}{4\hbar}}\Big)$.

Therefore, it remains to analyze the following integral
\begin{equation}\label{Oh}
\displaystyle\int_{y-a}^{y+a}e^{-\frac{(z-y)^2}{4\hbar}}e^{\frac{i\phi(z)}{\hbar}}dz.
\end{equation}

We define the operator  $L:=\frac{\hbar}{i}\frac{1}{\phi'(z)}\frac{d}{dz}$ and we note that $L^n(e^{\frac{i\phi}{\hbar}})=e^{\frac{i\phi}{\hbar}}$ for all $n\geq 1$. Then, integrating by parts we get

$$\int_{y-a}^{y+a}e^{-\frac{(z-y)^2}{4\hbar}}e^{\frac{i\phi(z)}{\hbar}}dz=$$
$$\int_{y-a}^{y+a}e^{-\frac{(z-y)^2}{4\hbar}}L(e^{\frac{i\phi(z)}{\hbar}})dz=\int_{y-a}^{y+a}\frac{\hbar}{i}\frac{e^{-\frac{(z-y)^2}{4\hbar}}}{\phi'(z)}\frac{d}{dz}(e^{\frac{i\phi(z)}{\hbar}})dz=$$

$$=\frac{\hbar}{i}\bigg( \frac{e^{-\frac{(z-y)^2}{4\hbar}}}{\phi'(z)}e^{\frac{i\phi(z)}{\hbar}}\bigg]_{y-a}^{y+a}-\int_{y-a}^{y+a}\frac{d}{dz}\bigg(\frac{e^{-\frac{(z-y)^2}{4\hbar}}}{\phi'(z)}\bigg)e^{\frac{i\phi(z)}{\hbar}}dz\bigg)= $$

$$=\frac{\hbar}{i}\bigg[ e^{-\frac{a^2}{4\hbar}} \bigg(\frac{e^{\frac{i\phi(y+a)}{\hbar}}}{\phi'(y+a)}-\frac{e^{\frac{i\phi(y-a)}{\hbar}}}{\phi'(y-a)}\bigg)$$
$$-\int_{y-a}^{y+a}\bigg(\frac{d}{dz} e^{-\frac{(z-y)^2}{4\hbar}}\frac{1}{\phi'(z)}-e^{-\frac{(z-y)^2}{4\hbar}}\frac{\phi''(z)}{\phi'(z)^2}\bigg)e^{\frac{i\phi(z)}{\hbar}}dz\bigg] $$

In order  to estimate the last integral, we will define the following functions $l^1_1(z)= \frac{1}{\phi'(z)}$ and $l^1_2(z)=\frac{\phi''(z)}{\phi'(z)^2}$. Note that these functions do not depend on $\hbar$, hence we can find a constant $C_1$, such that %
$|l^1_1(z)| \leq C_1$ and $|l^1_2(z)|\leq C_1$, if $z\in(y-a,y+a)$. Now, using the property $\Big|\int_I u(z) v(z)dz\Big|\leq \int_I |u(z)| |v(z)|dz \leq \max_{z\in I} |v(z)|\int_I |u(z)| dz$ and Lemma \ref{estderiv} we obtain

$$\bigg|\int_{y-a}^{y+a}\bigg(\frac{d}{dz} e^{-\frac{(z-y)^2}{4\hbar}}l_1^1(z)-e^{-\frac{(z-y)^2}{4\hbar}}l_2^1(z)\bigg)e^{\frac{i\phi(z)}{\hbar}}dz\bigg|\leq $$
$$\leq C_1 \int_{y-a}^{y+a}\bigg|\frac{d}{dz} e^{-\frac{(z-y)^2}{4\hbar}}\bigg|dz +C_1 \int_{y-a}^{y+a} e^{-\frac{(z-y)^2}{4\hbar}}dz  =O(1)$$


And finally,  the integral in \eqref{Oh} can be estimated by
$$\bigg|\int_{y-a}^{y+a}e^{-\frac{(z-y)^2}{4\hbar}}e^{\frac{i\phi(z)}{\hbar}}dz\bigg|\leq \hbar\bigg[\tilde C_1 e^{-\frac{a^2}{4\hbar}}+O(1)\bigg]=O(\hbar).$$

We want to show that the integral in \eqref{Oh} is $O(\hbar^N)$ for each $N\in\mathbb{N}$, and in order to obtain that we will integrate by parts several times.

Using the operator $L^2$ we get the following estimate
$$\int_{y-a}^{y+a}e^{-\frac{(z-y)^2}{4\hbar}}e^{\frac{i\phi(z)}{\hbar}}dz=$$
$$\int_{y-a}^{y+a}e^{-\frac{(z-y)^2}{4\hbar}}L^2(e^{\frac{i\phi(z)}{\hbar}})dz=\frac{\hbar^2}{i^2}\int_{y-a}^{y+a}\frac{e^{-\frac{(z-y)^2}{4\hbar}}}{\phi'(z)}\frac{d}{dz}\Big(\frac{1}{\phi'(z)}\frac{d}{dz}e^{\frac{i\phi(z)}{\hbar}}\Big)dz=$$
$$=\frac{\hbar^2}{i^2}\frac{e^{-\frac{(z-y)^2}{4\hbar}}}{\phi'(z)}\Big(\frac{1}{\phi'(z)}\frac{d}{dz}e^{\frac{i\phi(z)}{\hbar}}\Big)\bigg]_{y-a}^{y+a}-\frac{\hbar^2}{i^2}\int_{y-a}^{y+a}\frac{1}{\phi'(z)}\frac{d}{dz}\bigg(\frac{e^{-\frac{(z-y)^2}{4\hbar}}}{\phi'(z)}\bigg)\frac{d}{dz}e^{\frac{i\phi(z)}{\hbar}}dz=$$
\begin{equation}\label{h2}
 =\frac{\hbar}{i}\frac{e^{-\frac{(z-y)^2}{4\hbar}}}{\phi'(z)}e^{\frac{i\phi(z)}{\hbar}}\bigg]_{y-a}^{y+a}-\frac{\hbar^2}{i^2}\int_{y-a}^{y+a}\frac{1}{\phi'(z)}\frac{d}{dz}\bigg(\frac{e^{-\frac{(z-y)^2}{4\hbar}}}{\phi'(z)}\bigg)\frac{d}{dz}e^{\frac{i\phi(z)}{\hbar}}dz
\end{equation}
Now we will estimate the  integral in \eqref{h2}

$$\int_{y-a}^{y+a}\frac{1}{\phi'(z)}\frac{d}{dz}\bigg(\frac{e^{-\frac{(z-y)^2}{4\hbar}}}{\phi'(z)}\bigg)\frac{d}{dz}e^{\frac{i\phi(z)}{\hbar}}dz=$$
$$=\int_{y-a}^{y+a}\bigg(\frac{d}{dz} e^{-\frac{(z-y)^2}{4\hbar}}\frac{l_1^1(z)}{\phi'(z)}-e^{-\frac{(z-y)^2}{4\hbar}}\frac{l_2^1(z)}{\phi'(z)}\bigg)\frac{d}{dz}e^{\frac{i\phi(z)}{\hbar}}dz=$$
$$=\bigg(\frac{d}{dz} e^{-\frac{(z-y)^2}{4\hbar}}\frac{l_1^1(z)}{\phi'(z)}-e^{-\frac{(z-y)^2}{4\hbar}}\frac{l_2^1(z)}{\phi'(z)}\bigg)e^{\frac{i\phi(z)}{\hbar}}\bigg]_{y-a}^{y+a}- $$
\begin{equation}\label{i2}
-\int_{y-a}^{y+a}\frac{d}{dz}\bigg(\frac{d}{dz} e^{-\frac{(z-y)^2}{4\hbar}}\frac{l_1^1(z)}{\phi'(z)}-e^{-\frac{(z-y)^2}{4\hbar}}\frac{l_2^1(z)}{\phi'(z)}\bigg)e^{\frac{i\phi(z)}{\hbar}}dz
\end{equation}

and by calculating the derivative that appears  in the equation \eqref{i2} we have

$$\frac{d}{dz}\bigg(\frac{d}{dz} e^{-\frac{(z-y)^2}{4\hbar}}\frac{l_1^1(z)}{\phi'(z)}-e^{-\frac{(z-y)^2}{4\hbar}}\frac{l_2^1(z)}{\phi'(z)}\bigg) =$$
$$=\frac{d^2}{dz^2}e^{-\frac{(z-y)^2}{4\hbar}}\frac{l_1^1(z)}{\phi'(z)}+\frac{d}{dz}e^{-\frac{(z-y)^2}{4\hbar}}\frac{d}{dz}\frac{l_1^1(z)}{\phi'(z)}- \frac{d}{dz}e^{-\frac{(z-y)^2}{4\hbar}}\frac{l_2^1(z)}{\phi'(z)}-e^{-\frac{(z-y)^2}{4\hbar}}\frac{d}{dz}\frac{l_2^1(z)}{\phi'(z)}=$$
\begin{equation}\label{d2}
= \frac{d^2}{dz^2}e^{-\frac{(z-y)^2}{4\hbar}}l_1^2(z)+\frac{d}{dz}e^{-\frac{(z-y)^2}{4\hbar}}l^2_2-e^{-\frac{(z-y)^2}{4\hbar}}l^2_3,
\end{equation}
where $\displaystyle l_1^2=\frac{1}{\phi'^2},\,\, l_2^2=-\frac{3\phi''}{\phi'^3},\,\, l_3^2(z)=\frac{d}{dz}\frac{\phi''}{\phi'^3}=\frac{\phi'''\phi'-3\phi''^2}{\phi'^4}$ do not depend on $\hbar$ and are bounded by a constant $C_2$ in $(y-a,y+a)$ Lebesgue a.e.w.

Hence we can estimate the integral in \eqref{Oh} by substituting  the equation \eqref{d2} in \eqref{i2} and the equation \eqref{i2} in  \eqref{h2} in order to get

$$\bigg|\int_{y-a}^{y+a}e^{-\frac{(z-y)^2}{4\hbar}}e^{\frac{i\phi(z)}{\hbar}}dz\bigg|\leq $$
$$\leq \hbar[\tilde C_1 e^{-\frac{a^2}{4\hbar}}]+\hbar^2\bigg|\bigg(\frac{d}{dz} e^{-\frac{(z-y)^2}{4\hbar}}\frac{l_1^1(z)}{\phi'(z)}-e^{-\frac{(z-y)^2}{4\hbar}}\frac{l_2^1(z)}{\phi'(z)}\bigg)e^{\frac{i\phi(z)}{\hbar}}\bigg]_{y-a}^{y+a}\bigg|+ $$
$$
+\hbar^2C_2\bigg[\int_{y-a}^{y+a}\bigg|\frac{d^2}{dz^2}e^{-\frac{(z-y)^2}{4\hbar}}\bigg|dz+\int_{y-a}^{y+a}\bigg|\frac{d}{dz}e^{-\frac{(z-y)^2}{4\hbar}}\bigg|dz+\int_{y-a}^{y+a}e^{-\frac{(z-y)^2}{4\hbar}}dz\bigg]\leq $$
$$\leq \hbar[\tilde C_1 e^{-\frac{a^2}{4\hbar}}]+\hbar^2e^{-\frac{a^2}{4\hbar}}C\bigg(\frac{a}{\hbar}+1\bigg)+\hbar^2C_2[O(1)]=O(\hbar^2).$$
Note that in the last inequality we have used the Lemma \ref{estderiv}.

When we apply the operator $L^3$, we integrate by parts three times and hence  $\hbar^3$ appears multiplying all terms.
The only  term we need to be concern is
$$\int_{y-a}^{y+a}\frac{d}{dz}\bigg[\frac{1}{\phi'(z)}\bigg(\frac{d^2}{dz^2}e^{-\frac{(z-y)^2}{4\hbar}}l_1^2(z)+\frac{d}{dz}e^{-\frac{(z-y)^2}{4\hbar}}l^2_2-e^{-\frac{(z-y)^2}{4\hbar}}l^2_3\bigg)\bigg]e^{\frac{i\phi(z)}{\hbar}}dz,$$
 because all the other terms, as we see when we apply $L^2$, will have a factor $e^{-\frac{a^2}{4\hbar}}$ multiplying it.
Then, computing the derivative that appears in the  integrand we get

$$ \int_{y-a}^{y+a}\bigg[\frac{d^3}{dz^3}e^{-\frac{(z-y)^2}{4\hbar}}\frac{1}{{\phi'(z)^3}}-\frac{d^2}{dz^2}e^{-\frac{(z-y)^2}{4\hbar}}\frac{6\phi''}{\phi'^4}-\frac{d}{dz}e^{-\frac{(z-y)^2}{4\hbar}}\frac{4\phi'''\phi'-15\phi''^2}{\phi'^5}+$$
$$-e^{-\frac{(z-y)^2}{4\hbar}}\frac{\phi^{(4)}\phi'^2-10\phi'\phi''\phi'''+15\phi''^2}{\phi'^6}\bigg]
e^{\frac{i\phi(z)}{\hbar}}dz.$$

As  before,  we define $\displaystyle l_1^3=\frac{1}{{\phi'(z)^3}}$, $\displaystyle l_2^3=\frac{6\phi''}{\phi'^4}$, $\displaystyle l_3^3=\frac{4\phi'''\phi'-15\phi''^2}{\phi'^5}$, $\displaystyle l_4^3=\frac{\phi^{(4)}\phi'^2-10\phi'\phi''\phi'''+15\phi''^2}{\phi'^6}$ and we see that these functions are bounded by a constant $C_3$ that is independent of $\hbar$.

 Now using Lemma \ref{estderiv}, we obtain

$$\bigg|\int_{y-a}^{y+a}\frac{d}{dz}\bigg[\frac{1}{\phi'(z)}\bigg(\frac{d^2}{dz^2}e^{-\frac{(z-y)^2}{4\hbar}}l_1^2(z)+\frac{d}{dz}e^{-\frac{(z-y)^2}{4\hbar}}l^2_2-e^{-\frac{(z-y)^2}{4\hbar}}l^2_3\bigg)\bigg]e^{\frac{i\phi(z)}{\hbar}}dz\bigg|\leq$$
$$\leq C_3\bigg(\int_{y-a}^{y+a}\bigg|\frac{d^3}{dz^3}e^{-\frac{(z-y)^2}{4\hbar}}\bigg| dz+\int_{y-a}^{y+a}\bigg|\frac{d^2}{dz^2}e^{-\frac{(z-y)^2}{4\hbar}}\bigg| dz+$$
$$+\int_{y-a}^{y+a}\bigg|\frac{d}{dz}e^{-\frac{(z-y)^2}{4\hbar}}\bigg| dz+\int_{y-a}^{y+a}e^{-\frac{(z-y)^2}{4\hbar}} dz\bigg)=O\bigg(\frac{1}{\hbar}\bigg).$$
\bigskip

Therefore, because of the  multiplying factor $\hbar^3$, we get

$$\bigg|\int_{y-a}^{y+a}e^{-\frac{(z-y)^2}{4\hbar}}e^{\frac{i\phi(z)}{\hbar}}dz\bigg|=O(\hbar^2).$$

Analogously, when we apply $L^4$ it will appear the multiplicative term $\hbar^4$,  and now using Lemma \ref{estderiv} it can be shown that

$$\int_{y-a}^{y+a}\frac{d}{dz}\bigg(\frac{1}{\phi'(z)}\frac{d}{dz}\bigg(\frac{1}{\phi'(z)}\frac{d}{dz}\bigg(\frac{1}{\phi'(z)}\frac{d}{dz}\bigg(\frac{e^{-\frac{(z-y)^2}{4\hbar}}}{\phi'(z)}\bigg)\bigg)\bigg)\bigg)e^{\frac{i\phi(z)}{\hbar}}dz=O\bigg(\frac{1}{\hbar}\bigg),$$
hence
$$\bigg|\int_{y-a}^{y+a}e^{-\frac{(z-y)^2}{4\hbar}}e^{\frac{i\phi(z)}{\hbar}}dz\bigg|=O(\hbar^3).$$

Analogously, when we apply $L^5$ it will appear the multiplicative term $\hbar^5$, and now using Lemma \ref{estderiv} it can be shown that

$$\int_{y-a}^{y+a}\frac{d}{dz}\bigg(\frac{1}{\phi'(z)}\frac{d}{dz}\bigg(\frac{1}{\phi'(z)}\frac{d}{dz}
\bigg(\frac{1}{\phi'(z)}\frac{d}{dz}\bigg(\frac{1}{\phi'(z)}\frac{d}{dz}\bigg(\frac{e^{-\frac{(z-y)^2}{4\hbar}}}{\phi'(z)}
\bigg)\bigg)\bigg)\bigg)\bigg)e^{\frac{i\phi(z)}{\hbar}}dz=$$
$$O\bigg(\frac{1}{\hbar^2}\bigg),$$
hence
$$\bigg|\int_{y-a}^{y+a}e^{-\frac{(z-y)^2}{4\hbar}}e^{\frac{i\phi(z)}{\hbar}}dz\bigg|=O(\hbar^3).$$

Finally, it can be shown that if we apply $L^{2n}$ or  $L^{2n+1}$ we obtain
$$\bigg|\int_{y-a}^{y+a}e^{-\frac{(z-y)^2}{4\hbar}}e^{\frac{i\phi(z)}{\hbar}}dz\bigg|=O(\hbar^n),$$ and this implies
$$\bigg|\int_{y-a}^{y+a}e^{-\frac{(z-y)^2}{4\hbar}}e^{\frac{i\phi(z)}{\hbar}}dz\bigg|=O(\hbar^{\infty}).$$
\qed



{\bf Proof (of the Lemma \ref{estderiv}):}
 Let us first calculate the derivatives of the function  $e^{-\frac{(z-y)^2}{4\hbar}}$:
$$\frac{d}{dz}e^{-\frac{(z-y)^2}{4\hbar}}=e^{-\frac{(z-y)^2}{4\hbar}}\frac{(y-z)}{2\hbar},$$

$$\frac{d^2}{dz^2}e^{-\frac{(z-y)^2}{4\hbar}}=e^{-\frac{(z-y)^2}{4\hbar}}\bigg(\frac{(y-z)^2-2\hbar}{4\hbar^2}\bigg),$$

$$\frac{d^3}{dz^3}e^{-\frac{(z-y)^2}{4\hbar}}=e^{-\frac{(z-y)^2}{4\hbar}}\bigg(\frac{(y-z)^3-6\hbar(y-z)}{8\hbar^3}\bigg),$$

$$\frac{d^4}{dz^4}e^{-\frac{(z-y)^2}{4\hbar}}=e^{-\frac{(z-y)^2}{4\hbar}}\bigg(\frac{(y-z)^4-12\hbar(y-z)^2+12\hbar^2}{16\hbar^4}\bigg),$$

$$\frac{d^5}{dz^5}e^{-\frac{(z-y)^2}{4\hbar}}=e^{-\frac{(z-y)^2}{4\hbar}}\bigg(\frac{(y-z)^5-20\hbar(y-z)^3+60\hbar^2(y-z)}{32\hbar^5}\bigg),$$

$$\frac{d^6}{dz^6}e^{-\frac{(z-y)^2}{4\hbar}}=e^{-\frac{(z-y)^2}{4\hbar}}\bigg(\frac{(y-z)^6-30\hbar(y-z)^4+180\hbar^2(y-z)^2-120\hbar^3}{64\hbar^6}\bigg),$$
It is enough to see that  $\frac{d^{2n+1}}{dz^{2n+1}}e^{-\frac{(z-y)^2}{4\hbar}}|_{z=y}=0$ and  $\frac{d^{2n}}{dz^{2n}}e^{-\frac{(z-y)^2}{4\hbar}}|_{z=y}=\frac{C}{\hbar^n}$, where $C\neq 0$ is a constant.

Now we will estimate $\displaystyle\int_{y-a}^{y+a}\bigg|\frac{d^n}{dz^n} e^{-\frac{(z-y)^2}{4\hbar}}\bigg|dz,$ for $n\geq 1$.


Note that $\frac{d}{dz} e^{-\frac{(z-y)^2}{4\hbar}}$ only change the signal in $z=y$, then
$$\int_{y-a}^{y+a}\bigg|\frac{d}{dz} e^{-\frac{(z-y)^2}{4\hbar}}\bigg|dz = \int_{y-a}^{y}\frac{d}{dz} e^{-\frac{(z-y)^2}{4\hbar}}dz-\int_{y}^{y+a}\frac{d}{dz} e^{-\frac{(z-y)^2}{4\hbar}}dz = $$

$$=   e^{-\frac{(z-y)^2}{4\hbar}}\bigg|_{y-a}^{y}-e^{-\frac{(z-y)^2}{4\hbar}}\bigg|_{y}^{y+a}=2-2e^{-\frac{a^2}{4\hbar}}=O(1).$$

Observe that  $\frac{d^2}{dz^2}e^{-\frac{(z-y)^2}{4\hbar}}\Big|_{z=y}=-\frac 1 \hbar$ and it can  change of signal at most twice in $(y-a,y+a)$. Let us suppose that this happens, i.e., there exist $r_1<y< r_2$ points in $(y-a, y+a)$ such that the second derivative is positive in $(y-a, r_1)$ and in $(r_2, y+a)$ and negative in $(r_1, r_2)$, then we obtain

$$\int_{y-a}^{y+a}\bigg|\frac{d^2}{dz^2}e^{-\frac{(z-y)^2}{4\hbar}}\bigg|dz=$$
$$=\int_{y-a}^{r_1}\frac{d^2}{dz^2}e^{-\frac{(z-y)^2}{4\hbar}}dz-\int_{r_1}^{r_2}\frac{d^2}{dz^2}e^{-\frac{(z-y)^2}{4\hbar}}dz+\int_{r_2}^{y+a}\frac{d^2}{dz^2}e^{-\frac{(z-y)^2}{4\hbar}}dz=$$
$$=\frac{d}{dz}e^{-\frac{(z-y)^2}{4\hbar}}\bigg|_{y-a}^{r_1}-\frac{d}{dz}e^{-\frac{(z-y)^2}{4\hbar}}\bigg|_{r_1}^{r_2}+\frac{d}{dz}e^{-\frac{(z-y)^2}{4\hbar}}\bigg|_{r_2}^{y+a}=$$
$$=2\bigg(e^{-\frac{(r_1-y)^2}{4\hbar}}\bigg(\frac{y-r_1}{\hbar}\bigg)-e^{-\frac{a^2}{4\hbar}}\bigg(\frac{a}{\hbar}\bigg)-e^{-\frac{(r_2-y)^2}{4\hbar}}\bigg(\frac{y-r_2}{\hbar}\bigg)\bigg),$$
 and this decay fast to zero as $\hbar\to 0$.

Note that  $\frac{d^3}{dz^3}e^{-\frac{(z-y)^2}{4\hbar}}\Big|_{z=y}=0$ and it can  change of signal at most three times in $(y-a,y+a)$. Let us suppose that this happens, i.e., there exist $r_1<y< r_2$ points in $(y-a, y+a)$ such that the  derivative of order 3 is positive in $(y-a, r_1)$ and in $(y,r_2)$  and negative in $(r_1, y)$ and in $(r_2, y+a)$, then, we obtain
$$\int_{y-a}^{y+a}\bigg|\frac{d^3}{dz^3}e^{-\frac{(z-y)^2}{4\hbar}}\bigg|dz=$$
$$=\int_{y-a}^{r_1}\frac{d^3}{dz^3}e^{-\frac{(z-y)^2}{4\hbar}}dz-\int_{r_1}^{y}\frac{d^3}{dz^3}e^{-\frac{(z-y)^2}{4\hbar}}dz+$$
$$+\int_{y}^{r_2}\frac{d^3}{dz^3}e^{-\frac{(z-y)^2}{4\hbar}}dz-\int_{r_2}^{y+a}\frac{d^3}{dz^3}e^{-\frac{(z-y)^2}{4\hbar}}dz=$$
$$=\frac{d^2}{dz^2}e^{-\frac{(z-y)^2}{4\hbar}}\bigg|_{y-a}^{r_1}-\frac{d^2}{dz^2}e^{-\frac{(z-y)^2}{4\hbar}}\bigg|_{r_1}^{y}+\frac{d^2}{dz^2}e^{-\frac{(z-y)^2}{4\hbar}}\bigg|_{y}^{r_2}-\frac{d^2}{dz^2}e^{-\frac{(z-y)^2}{4\hbar}}\bigg|_{r_2}^{y+a}=$$
$$=2e^{-\frac{(r_1-y)^2}{4\hbar}}\bigg(\frac{(y-r_1)^2-\hbar}{\hbar^2}\bigg)-2e^{-\frac{a^2}{4\hbar}}\bigg(\frac{a^2-\hbar}{\hbar^2}\bigg)+$$
$$+\frac{2}{\hbar}+2e^{-\frac{(r_2-y)^2}{4\hbar}}\bigg(\frac{(y-r_2)^2-\hbar}{\hbar^2}\bigg)=O\bigg(\frac{1}{\hbar}\bigg).$$
It became clear from these calculations that,  $\displaystyle\int_{y-a}^{y+a}\bigg|\frac{d^{2n}}{dz^{2n}}e^{-\frac{(z-y)^2}{4\hbar}}\bigg|dz$ decay fast to zero as $\hbar\to 0$, because $\displaystyle\frac{d^{2n}}{dz^{2n}}e^{-\frac{(z-y)^2}{4\hbar}}$ does not change of signal in $z=y$. Moreover,  the only term in  $\displaystyle\int_{y-a}^{y+a}\bigg|\frac{d^{2n+1}}{dz^{2n+1}}e^{-\frac{(z-y)^2}{4\hbar}}\bigg|dz$  that does not decay fast to zero is $\frac{d^{2n}}{dz^{2n}}e^{-\frac{(z-y)^2}{4\hbar}}\bigg|_{z=y}=\frac{C}{\hbar^{n}}, $ and this implies  $$\int_{y-a}^{y+a}\bigg|\frac{d^{2n+1}}{dz^{2n+1}}e^{-\frac{(z-y)^2}{4\hbar}}\bigg|dz=O\bigg(\frac{1}{\hbar^n}\bigg).$$
\qed

\vspace{0.5cm}

\newpage

{\bf Instituto de Matematica - UFRGS - Brasil}
\medskip

A. O. Lopes was partially supported by CNPq and INCT and J. Mohr was partially supported by CNPq.

\end{document}